# A general learning algorithm for functions between metric spaces

Kerry M. Soileau

August 16, 2007

ABSTRACT

In this paper we show how to approximate ("learn") a function $f : X \to Y$, where $\langle X, \sigma \rangle$ and $\langle Y, \rho \rangle$ are metric spaces.

## 1. INTRODUCTION

Given a function $f : X \to Y$, where $\langle X, \sigma \rangle$ and $\langle Y, \rho \rangle$ are metric spaces, we seek to approximate $f$ through a sequence of evaluations. For every nonempty finite $A \subseteq X$, define $\wp(x, A) \equiv \{ t \in A; \forall u \in A, \sigma(t, x) \leq \sigma(u, x) \}$, and let

$u(S) \equiv$ a uniformly randomly chosen member of a finite set $S$. Let $\{x_n\}_{n=1}^{\infty} \subseteq X$ be a sequence of distinct elements, and let $\{f(x_n)\}_{n=1}^{\infty} \subseteq Y$ be its image. Fix numbers $\varepsilon \in (0, \infty)$ and $q \in [\frac{1}{2}, 1)$. Define $A_n \subseteq X$ according to $A_1 = \{x_1\}$ and

$$A_n = \begin{cases} A_{n-1} \cup \{x_n\} & \text{if } \rho\big(f\big(u\big(\wp(x_n, A_{n-1})\big)\big), f(x_n)\big) > \varepsilon \\ A_{n-1} \sim \{u(\wp(x_n, A_{n-1}))\} & \text{with probability } \dfrac{1}{q} - 1 \text{ if } \rho\big(f\big(u\big(\wp(x_n, A_{n-1})\big)\big), f(x_n)\big) \leq \varepsilon \\ A_{n-1} & \text{with probability } 2 - \dfrac{1}{q} \text{ if } \rho\big(f\big(u\big(\wp(x_n, A_{n-1})\big)\big), f(x_n)\big) \leq \varepsilon \end{cases}$$

for $n > 1$.

## 2. THE APPROXIMATION FUNCTION

We now define a sequence of functions $\{f_n : X \to Y\}_{n=1}^{\infty}$. Define $f_n(x) \equiv f\big(u\big(\wp(x, A_n)\big)\big)$ for $n = 1, 2, 3, \cdots$

A simple calculation shows the growth behavior of $A_n$.

$$|A_n|-|A_{n-1}| = \begin{cases} 1 & \text{if } \rho(f_{n-1}(x_n), f(x_n)) > \varepsilon \\ -1 & \text{with probability } \dfrac{1}{q}-1 \text{ if } \rho(f_{n-1}(x_n), f(x_n)) \le \varepsilon \\ 0 & \text{with probability } 2-\dfrac{1}{q} \text{ if } \rho(f_{n-1}(x_n), f(x_n)) \le \varepsilon \end{cases}$$

So

$$E\big(|A_n|-|A_{n-1}| \big| \rho(f_{n-1}(x_n), f(x_n)) > \varepsilon \big) = 1$$

$$E\big(|A_n|-|A_{n-1}| \big| \rho(f_{n-1}(x_n), f(x_n)) \le \varepsilon \big) = 1 - \frac{1}{q}$$

Hence

$$E(|A_n|-|A_{n-1}|) = 1 \cdot P\big(\rho(f_{n-1}(x_n), f(x_n)) > \varepsilon\big) + \left(1-\frac{1}{q}\right) P\big(\rho(f_{n-1}(x_n), f(x_n)) \le \varepsilon\big)$$

$$= 1 - P\big(\rho(f_{n-1}(x_n), f(x_n)) \le \varepsilon\big) + \left(1-\frac{1}{q}\right) P\big(\rho(f_{n-1}(x_n), f(x_n)) \le \varepsilon\big)$$

$$= 1 - \frac{1}{q} P\big(\rho(f_{n-1}(x_n), f(x_n)) \le \varepsilon\big)$$

<u>Theorem</u>: If $\lim_{n \to \infty} E(|A_n|-|A_{n-1}|) = 0$, then $\lim_{n \to \infty} P\big(\rho(f_{n-1}(x_n), f(x_n)) \le \varepsilon\big) = q$.

<u>Proof</u>:

$$0 = \lim_{n \to \infty} E(|A_n|-|A_{n-1}|) = \lim_{n \to \infty} \left(1 - \frac{1}{q} P\big(\rho(f_{n-1}(x_n), f(x_n)) \le \varepsilon\big)\right)$$

$$= 1 - \frac{1}{q} \lim_{n \to \infty} P\big(\rho(f_{n-1}(x_n), f(x_n)) \le \varepsilon\big)$$

Hence $\lim_{n \to \infty} P\big(\rho(f_{n-1}(x_n), f(x_n)) \le \varepsilon\big) = q$.

This suggests that if we attain convergence of the expected size of $A_n$ in the limit, i.e. $\lim_{n\to\infty} E(|A_n| - |A_{n-1}|) = 0$, we must expect to get $\lim_{n\to\infty} P(\rho(f_{n-1}(x_n), f(x_n)) \leq \varepsilon) = q$.

Therefore we may interpret the parameter $q$ as the limiting "performance" of the approximation, i.e. the limiting probability that the approximation will return a value within distance $\varepsilon$ of the correct value.

## 3. REFERENCES

[1]     E.T. Whittaker and G.N. Watson, *A Course of Modern Analysis*, Cambridge University Press, London and New York, 1958.

International Space Station Program Office, Avionics and Software Office, Mail Code OD, NASA Johnson Space Center
Houston, TX 77058
E-mail address: ksoileau@yahoo.com